\documentclass[11pt,twoside]{article}

\usepackage{a4wide}
\usepackage{amsfonts}
\usepackage{amssymb}
\usepackage{amsmath}
\usepackage{graphicx}
\usepackage{nopageno}

\newcommand{\ignore}[1]{}

\def\@begintheorem#1#2{\par\bgroup{\sc #1\ #2. }\it\ignorespaces}
\def\@opargbegintheorem#1#2#3{\par\bgroup{\sc #1\ #2\ (#3). } \it\ignorespaces}
\def\@endtheorem{\egroup}
\newtheorem{theorem}{Theorem}[section]
\newtheorem{corollary}[theorem]{Corollary}
\newtheorem{lemma}[theorem]{Lemma}

\newtheorem{example}[theorem]{Example}
\newtheorem{proposition}[theorem]{Proposition}
\newtheorem{definition}[theorem]{Definition}
\newcommand{\bt}[1]{\begin{theorem}\label{#1}}
\newcommand{\bc}[1]{\begin{corollary}\label{#1}}
\newcommand{\bl}[1]{\begin{lemma}\label{#1}}
\newcommand{\be}[1]{\begin{example}\label{#1}}
\newcommand{\bp}[1]{\begin{proposition}\label{#1}}
\newcommand{\ba}[1]{\begin{algorithm}\rm\label{#1}}
\newcommand{\bd}[1]{\begin{definition}\rm\label{#1}}{\normalsize }
\newcommand{\bpr}{\noindent {\em Proof. }}
\newcommand{\et}{\end{theorem}}
\newcommand{\ec}{\end{corollary}}
\newcommand{\el}{\end{lemma}}
\newcommand{\ee}{\end{example}}
\newcommand{\ep}{\end{proposition}}
\newcommand{\ed}{\end{definition}}
\newcommand{\epr}{{\ \vbox{\hrule\hbox{%
\vrule height1.3ex\hskip0.8ex\vrule}\hrule}}\\\par}

\def\Z{\mathbb{Z}}
\def\nk{\{0,1\}^n_k}
\def\nt{\{0,1\}^n_3}
\def\sign{{\rm sign}}

\begin{document}

\title{\bf Hypergraphic Degree Sequences are Hard}

\author{Antoine Deza, Asaf Levin, Syed M. Meesum, and Shmuel Onn}
\date{}

\maketitle

\noindent
A {\em $k$-hypergraph on $[n]$} is a subset $H\subseteq\nk:=\{x\in\{0,1\}^n:\|x\|_1=k\}$.
The {\em degree sequence of $H$} is the vector $d=\sum H:=\sum\{x:x\in H\}$.
We consider the following decision problem: given $k$ and $d\in\Z_+^n$,
is $d$ the degree sequence of some hypergraph $H\subseteq\nk$?
For $k=2$ (graphs) the celebrated work of Erd\H{o}s and Gallai \cite[1960]{EG}
implies that $d$ is a degree sequence of a graph if and only if $\sum d_i$ is even and,
permuting $d$ so that $d_1\geq\dots\geq d_n$, the inequalities
$\sum_{i=1}^j d_i-\sum_{i=l+1}^n d_i\leq j(l-1)$ hold for $1\leq j\leq l\leq n$,
yielding a polynomial time algorithm. For $k=3$ the problem was raised over 30 years ago
by Colbourn, Kocay, and Stinson \cite[1986, Problem 3.1]{CKS} and was solved by
Deza, Levin, Meesum, and Onn \cite[2018]{DLMO}. 

Here is the statement and its short proof.

\vskip.2cm
\noindent{\bf Theorem:}
It is NP-complete to decide if $d\in\Z_+^n$ is the degree sequence of a $3$-hypergraph.

\vskip.2cm
\bpr
The problem is in NP since if $d$ is a degree sequences then a hypergraph
$H\subseteq\nt$ of cardinality $|H|\leq{n\choose3}=O(n^3)$
can be exhibited and $d=\sum H$ verified in polynomial time.

Let ${\bf 1}$ be the all-ones vector. We consider the following three decision problems.

\vskip.2cm
\noindent
(1) Given $a\in\Z_+^n$, $b\in\Z_+$ with $3{\bf 1}a=nb$, is there
$F\subseteq\{x\in\nt:ax=b\}$ with $\sum F={\bf 1}$ ?

\noindent
(2) Given $w\in\Z^n$, $c\in\Z_+^n$ with $wc=0$,
is there $G\subseteq\{x\in\nt:wx=0\}$ with $\sum G=c$ ?

\noindent
(3) Given $d\in\Z_+^n$, is there $H\subseteq\nt$ with $\sum H=d$ ?

\vskip.2cm
Problem (1) is the so-called {\em $3$-partition problem} which is known to be NP-complete \cite{GJ}.

First we reduce (1) to (2). Given $a,b$ with $3{\bf 1}a=nb$, let $w:=3a-b{\bf 1}$ and
$c:={\bf 1}$. Then $wc=0$. Now, for any $x\in\nt$ we have $wx=3ax-b{\bf 1}x=3(ax-b)$
so $x$ satisfies $ax=b$ if and only if $wx=0$.
So the answer to (1) is YES if and only if the answer to (2) is YES.

Second we reduce (2) to (3). Given $w,c$, with $wc=0$, define $d:=c+\sum S_+$, where
$S_\sigma:=\{x\in\nt:\sign(wx)=\sigma\}$ for $\sigma=-,0,+$.
Suppose there is a $G\subseteq S_0$ with $\sum G=c$. Then $H:=G\cup S_+$ satisfies $\sum H=d$.
Suppose there is an $H\subseteq\nt$ with $\sum H=d$. Then
$$w\sum S_+\ =\ w(c+\sum S_+)\ =\ w\sum H\ =\
w\sum(H\cap S_-)+w\sum(H\cap S_0)+w\sum(H\cap S_+)$$
which implies $H\cap S_-=\emptyset$ and $H\cap S_+=S_+$.
Therefore $G:=H\cap S_0$ satisfies $\sum G=\sum H-\sum S_+=c$.
So the answer to (2) is YES if and only if the answer to (3) is YES.
\epr

\end{document}